\documentclass{siamltex}
\usepackage[latin1]{inputenc}

\usepackage{pdfpages} 

\usepackage[T1]{fontenc} 

\usepackage{fancyhdr}
\fancypagestyle{empty}{}%
\usepackage{graphicx,xcolor}
\usepackage{subfigure}          

\usepackage[english]{babel}
\usepackage[autostyle=true,german=quotes]{csquotes}

\usepackage{epigraph}

\usepackage{placeins}

\usepackage{tabularx}

\usepackage{amsmath}
\usepackage{amsfonts}
\usepackage{amssymb}
\usepackage{mathabx}

\usepackage{lmodern,textcomp}

\usepackage{algorithm}
\usepackage{algpseudocode}

\usepackage{enumitem}
\newlist{Aufz}{enumerate}{10}
\setlist[Aufz]{label*=\arabic*.}

\usepackage{xspace}
\newcommand\largeparbreak{\par\bigskip}

\newcommand{\mdtheorem}[2]{\newtheorem{#1}{#2}}
\mdtheorem{Definition}{Definition}
\mdtheorem{Remark}{Remark}
\mdtheorem{Theorem}{Theorem}
\mdtheorem{Conjecture}{Conjecture}
\mdtheorem{Lemma}{Lemma}
\mdtheorem{Example}{Example}

\newcommand{\R}{\mathbb{R}\xspace}
\newcommand{\N}{\mathbb{N}\xspace}
\newcommand{\C}{\mathbb{C}\xspace}
\newcommand{\Cno}{\mathbb{C}\setminus \lbrace 0 \rbrace \xspace}

\newcommand{\cM}{\mathcal{M}\xspace}

\newcommand{\cP}{\mathcal{P}\xspace}
\newcommand{\cQ}{\mathcal{Q}\xspace}
\newcommand{\cG}{\mathcal{G}\xspace}

\newcommand{\tcC}{\tilde{\mathcal{C}}\xspace}
\newcommand{\cU}{\mathcal{U}\xspace}

\newcommand{\cK}{\mathcal{K}\xspace}

\newcommand{\cT}{\mathcal{T}\xspace}

\newcommand{\IDR}{\textsf{IDR($s$)}\xspace}
\newcommand{\IDReins}{\textsf{IDR($1$)}\xspace}
\newcommand{\IDRstab}{\textsf{IDR($s$)stab($\ell$)}\xspace}

\newcommand{\IDRO}{\textsf{IDR}\xspace}

\newcommand{\SRIDR}{\textsf{SRIDR($s$)}\xspace}

\newcommand{\Mstab}{\textsf{$\mathcal{M}$($s$)stab($\ell$)}\xspace}

\newcommand{\BiCG}{\textsf{BiCG}\xspace}
\newcommand{\BiCGstab}{\textsf{BiCGStab}\xspace}
\newcommand{\BiCGstabL}{\textsf{BiCGStab($\ell$)}\xspace}

\newcommand{\GMRES}{\textsf{GMRES}\xspace}

\newcommand{\GCR}{\textsf{GCR}\xspace}
\newcommand{\GCRO}{\textsf{GCRO}\xspace}
\newcommand{\MINRES}{\textsf{MINRES}\xspace}
\newcommand{\RMINRES}{\textsf{R-MINRES}\xspace}

\newcommand{\GCRODR}{\textsf{GCRO-DR}\xspace}
\newcommand{\GCROT}{\textsf{GCR-OT}\xspace}

\newcommand{\MV}{\textsf{MV}\xspace}

\newcommand{\hMV}{\textsf{\#MVs}\xspace}
\newcommand{\hRD}{\textsf{\#RDs}\xspace}

\newcommand{\rhs}{\textsf{RHS}\xspace}

\newcommand{\bA}{\textbf{A}\xspace}

\newcommand{\bU}{\textbf{U}\xspace}
\newcommand{\bV}{\textbf{V}\xspace}
\newcommand{\bZ}{\textbf{Z}\xspace}

\newcommand{\bI}{\textbf{I}\xspace}
\newcommand{\bR}{\textbf{R}\xspace}
\newcommand{\bL}{\textbf{L}\xspace}

\newcommand{\bz}{\textbf{z}\xspace}
\newcommand{\be}{\textbf{e}\xspace}

\newcommand{\bx}{\textbf{x}\xspace}

\newcommand{\by}{\textbf{y}\xspace}
\newcommand{\bb}{\textbf{b}\xspace}
\newcommand{\br}{\textbf{r}\xspace}

\newcommand{\bg}{\textbf{g}\xspace}
\newcommand{\bt}{\textbf{t}\xspace}

\newcommand{\bq}{\textbf{\textit{q}}\xspace}

\newcommand{\bu}{\textbf{u}\xspace}
\newcommand{\bv}{\textbf{v}\xspace}

\newcommand{\bp}{\textbf{p}\xspace}

\newcommand{\tbx}{\tilde{\textbf{x}}\xspace}

\newcommand{\tbr}{\tilde{\textbf{r}}\xspace}
\newcommand{\tbv}{\tilde{\textbf{v}}\xspace}

\newcommand{\s}{$s$}

\newcommand{\bO}{\textbf{0}\xspace}

\newcommand{\nEqns}{{n_{\text{Systems}}}\xspace}

\newcommand{\cond}{\operatorname{cond}\xspace}
\newcommand{\opspan}{\operatorname{span}\xspace}
\newcommand{\opImage}{\operatorname{range}\xspace}

\newcommand{\ha}{^{(1)}\xspace}
\newcommand{\hb}{^{(2)}\xspace}
\newcommand{\hia}{^{(\iota)}\xspace}
\newcommand{\hib}{^{(\iota+1)}\xspace}

\newcommand{\tol}{\mathrm{tol}\xspace}

\newcommand{\I}{we\xspace}
\newcommand{\gI}{We\xspace}
\newcommand{\my}{our\xspace}

\newcommand{\We}{We\xspace}
\newcommand{\we}{we\xspace}

	{\end{cases}\right.}

\title{M(\s )stab($\ell$): A Generalization of IDR(\s )stab($\ell$) for Sequences of Linear Systems}
\author{Martin P. Neuenhofen\thanks{RWTH Aachen University, Germany (\texttt{Martin.Peter.Neuenhofen@rwth-aachen.de}).}
	}

\date{\today}
\begin{document}


\maketitle

\begin{abstract}
\gI propose \Mstab, a novel Krylov subspace recycling method for the iterative solution of sequences of linear systems with fixed system matrix and changing right-hand sides. This new method is a straight and simple generalization of \IDRstab. \IDRstab in turn is a very efficient method and generalization of \BiCGstab.

The theory of \Mstab is based on a generalization of the \IDRO theorem and Sonneveld spaces.

Numerical experiments indicate that \Mstab can solve sequences of linear systems faster than its corresponding \IDRstab variant. Instead, when solving a single system both methods are identical.
\end{abstract}

\begin{AMS}
Primary, 
93C05, 
65F10; 
Secondary,
93A15, 
65F50, 
65N22, 
76M10. 
\end{AMS}

\begin{keywords}
  Sequence of linear systems, iterative solvers, Krylov subspaces, recycling, short recurrences, IDRstab, BiCGStab, SRIDR, Sonneveld spaces.
\end{keywords}


\section{Introduction}
\We consider iterative methods for the solution of sequences of large sparse symmetric and nonsymmetric linear systems
\begin{align}
	\bA \cdot \bx\hia = \bb\hia,\quad \iota = 1,...,\nEqns\,, \label{eqn:AxbSequence}
\end{align}
with fixed regular $\bA \in \C^{N \times N}$, where the right-hand sides (\rhs) $\bb\hib \in \C^N$  depend on the former solution $\bx\hia$ to the \rhs $\bb\hia$. Thus the systems must be solved \textit{one after the other}. Such situations occur e.g. when applying an implicit time stepping scheme to numerically solve a non-stationary partial differential equation. Areas of application are e.g. topology optimization \cite{RMINRES}, model reduction \cite{RGCRO}, structural dynamics \cite{GCRO-DR}, circuit analysis \cite{Compressing} and fluid dynamics \cite{NavStokes}. In all these referenced works a technique called {Krylov subspace recycling} is used.

\subsection{Krylov Subspace Recycling}
The idea of Krylov Subspace Recycling is to keep information from a former solution process to solve a subsequent system more efficiently. Imagine to solve a system $\bA \cdot \bx = \bb\ha$ using a Krylov subspace method with, e.g., $100$ iterations. During the solution process, a basis matrix for a $100$-dimensional search space $\cU$ (which is a Krylov subspace) is built.

When afterwards solving a subsequent system $\bA \cdot \bx = \bb\hb$, information from the old search space $\cU$ can be reused to possibly reduce computational effort for the solution to $\bb\hb$.

\subsubsection{Literature Review}
Most Krylov subspace recycling methods are based on the concept of augmenting the recycled search space $\cU$ iteratively: With a \emph{recycling space} $\cU$ from the solution of a former system, the current system is solved by adding new directions to $\cU$, obtaining a larger space $\hat{\cU}$. Then these methods compute a projection solution to the current system in this augmented search space $\hat{\cU}$, see e.g. \cite{RGCRO}.

The size of the recycling space $\cU$ is limited because of storage requirements for its basis and growing computational effort due to long recurrences in the computation of augmented directions.

Comparable to restarted \GMRES, in \textsf{GMRES-DR} \cite{RGMRES} this problem is treated by a warm restart: whenever the search space $\hat{\cU}$ becomes too large, a smaller subspace $\cU$ is extracted from it by deflation with Ritz vectors. The space $\cU$ is then used as recycling space for a restart and iteratively augmented. This can be done for one single system as well as for a sequence of systems. For the case of Hermitian systems the nested \GMRES precedure is replaced by \MINRES, yielding the recycling method \RMINRES \cite{RMINRES}. 

There also exist recycling variants for \GCR based on similar concepts. E.g., \GCRODR \cite{GCRO-DR} is a \GCRO method with a recycled search space $\cU$ obtained by deflation with Ritz vectors. \GCROT \cite{GCR-OT} chooses the recycled space $\cU$ by optimal truncation, instead. 

Besides there exist recycling variants of \BiCG and \BiCGstab \cite{RBiCG,R-BiCG3,RBiCGstab}. In contrast to the above long-recurrence methods, for \BiCG and \BiCGstab a projection basis must additionally be stored. This is due to the reason that the test space of the nonsymmetric Lanczos process differs from the search space, whereas in the Arnoldi process both spaces are similar.

\subsubsection{Outline}
Whereas conventional Krylov subspace methods iteratively build successively {growing} search spaces in which they compute a solution, Induced-Dimension-Reduction (\IDRO) methods iteratively restrict the residual of a numerical solution into successively {shrinking} subspaces.

Whereas conventional Krylov subspace recycling methods \textit{augment} their search space such that it is {greater} right from the beginning, \I propose in this paper an \IDRO-method that instead directly starts with a residual that lies in a {smaller} subspace. Thus \my recycling concept is the antipode to the augmentation approach.

\subsection{New Contributions}
In this paper \I present a novel short-recurrence iterative method that is based on recycling of shrinking nested subspaces. This new method ensures that the residual does not grow during the recycling phase. It can easily be incorporated into existent implementations of all \IDRO-methods. In constrast to the recycling methods described above, it can recycle large test spaces with only a limited number of stored basis vectors. Finally \my new method is a generalization of \IDRstab, the most efficient iterative short-recurrence method for single nonsymmetric linear systems known so far \cite[l.1-2]{IDRstab-Paper}. Side contributions are given on properties of Sonneveld spaces and their orthogonal complements.

\subsection{Structure}
Section 2 reviews graphically the foundations of \IDRO (with side contributions) and derives \IDRstab from scratch. Then in section 3 \I motivate \my new method by reviewing and discussing two early recycling approaches from literature, specifically the method of Miltenberger \cite{IDR-Milten} and \my \SRIDR \cite[sec. 2]{Report15}. Finally, in section 4 \I generalize the \IDRO theorem, build on that the new method \Mstab, and compare it to the earlier recycling methods from section 3. In section 5 numerical experiments demonstrate the efficiency and fast termination behaviour of \Mstab.

\section{Induced Dimension Reduction (\IDRO)}
This section gives both an introduction to \IDRO and provides new insight into the theory of Sonneveld spaces.

\subsection{Geometric Idea}
\enquote{\BiCGstab is the most popular short-recurrence iterative method for solving [single] large nonsymmetric systems of equations} \cite[Sec.1, l.1]{IDRstab-Paper}. It is algebraically equivalent to the iterative method \IDReins \cite{IDReins}. The latter method provides a simple termination theory based on {Sonneveld spaces}. To discuss these spaces \I use the following definitions.
\begin{Definition}[Krylov Subspace \& Block Krylov Subspace]
	Given $\bA \in \C^{N \times N}$ and $\bt \in \C^N$, the \emph{Krylov subspace} of \emph{level $j$} is
	\begin{align*}
		\cK_j(\bA;\bt) := \operatornamewithlimits{span}_{k=0,...,j-1}\lbrace \bA^k \cdot \bt\rbrace\,.
	\end{align*}
	$\bt$ is called \emph{starting vector}. This is the conventional definition. However, \I will only use \emph{block Krylov subspaces}. In these the starting vector is replaced by a vector space $\cT < \C^N$, i.e.
	\begin{align*}
	\cK_j(\bA;\cT) := \operatornamewithlimits{span}_{k=0,...,j-1}\lbrace \bA^k \cdot \cT\rbrace\,.
	\end{align*}
\end{Definition}
\begin{Definition}[Sonneveld Spaces]
	Given $\bA \in \C^{N \times N}$, a vector space $\cP < \C^N$ and a sequence $\lbrace \omega_j \rbrace_{j \in \N} \subset \Cno$, define recursively $\cG_0 := \C^N$ and
	\begin{align}
		\cG_{j} := (\bI - \omega_j \cdot \bA) \cdot (\cG_{j-1} \cap \cP^\perp)\quad \forall \, j \in \N\,.\label{eqn:DefSonneveld}
	\end{align}
	The vector spaces $\cG_j$ are called \emph{Sonneveld spaces} of \emph{level} $j$, respectively. \gI call the $\omega_j$ \emph{relaxations} and $\cP$ the \emph{cut-space}\footnote{Originally $\cP$ was called \emph{shadow (residual) space}.}.
\end{Definition}
\begin{Definition}[Stabilization Polynomial]
	Given the relaxations $\lbrace \omega_j \rbrace_{j \in \N}$ from above, define the \emph{stabilization polynomials}
	\begin{align*}
		p_{i,j}(t) := \prod_{k=i+1}^j (1 - \omega_k \cdot t)\quad\forall\ j>i \in \N_0\,.
	\end{align*}
\end{Definition}

\begin{figure}
	\centering
	\includegraphics[width=0.8\linewidth]{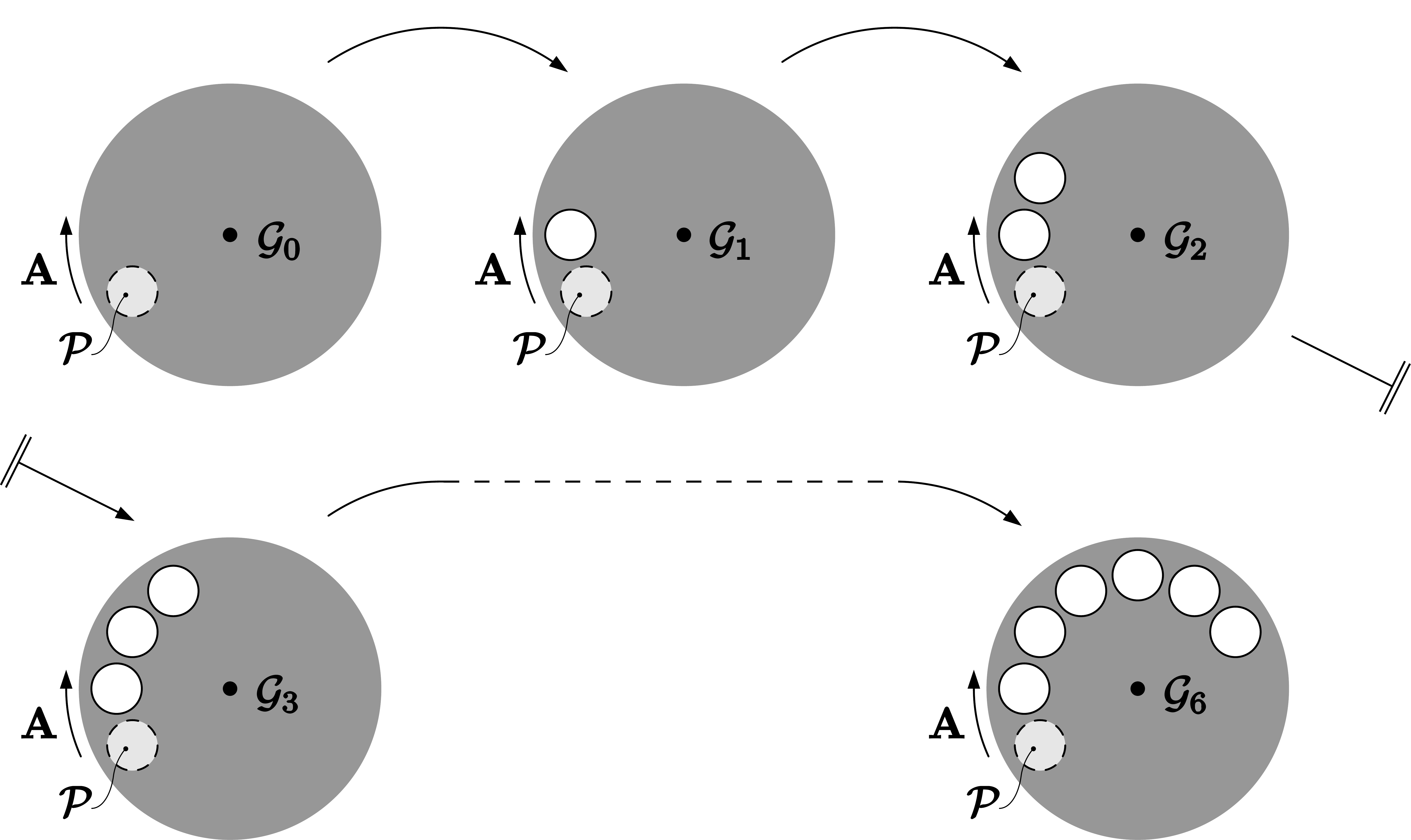}
	\caption{\textit{Geometric intution of induced dimension reduction by rotary dials.}}
	\label{fig:DialPlate}
\end{figure}

The geometric principle of \IDRO-methods is to move a residual $\br$ of a numerical solution iteratively into higher level Sonneveld spaces. With a higher level the Sonneveld space becomes smaller and eventually collapses to $\lbrace \bO \rbrace$. This leads to the exact solution due to $\br=\bO$ (details in Section 2.3). All Krylov subspace methods aim for such finite termination properties as these cause the favored superlinear convergence behaviour.
\largeparbreak
\paragraph{Intuition}
\We can think of Sonneveld spaces as rotary dials such as shown in Fig.~\ref{fig:DialPlate}. Figuratively speaking, when producing a rotary dial \we take a round plate (comparable to $\cG_0$, top left) and drill a borehole into it at a specific position, e.g. where dial number ones lies. Mathematically the bore position and shape is defined by the cut-space $\cP$, light grey in the figure. The drill process is mathematically comparable to cutting $\cP$ out of the Sonneveld space $\cG_0$. Figuratively, afterwards the plate is rotated such that the borehole moves from dial number one to dial number two. Mathematically, the rotation is comparable to multiplying\footnote{Let us neglect the relaxations $\omega_j$ for the moment.} the cutted Sonneveld space with $\bA$. To draw a bridge from our visual object to mathematics, the drilled and rotated dial plate equates $\cG_1$. Repeating this procedure of drilling and rotating, the dial plates resp. Sonneveld spaces look as in Fig.~\ref{fig:DialPlate}.

\We see from the picture that the material domain of the dial plate after a subsequent production step is contained in the material domain of each former plate. The same also holds for the Sonneveld spaces, i.e. $\cG_j \subseteq \cG_{j-1} \subseteq ... \subseteq \cG_0$. Additionally \we see that by drilling holes into the plate it loses a certain amount of material. Analogously the dimensions of the Sonneveld spaces decrease in the canonical case by a certain amount, i.e. $\dim(\cG_j) \leq \dim(\cG_{j-1}) - \dim(\cP)$ for $j=1,2,...$ until the space vanishes. The \IDRO theorem \cite{IDReins,IDR-report,IDR-Gutkn} proofs both mentioned properties of the Sonneveld spaces, and later \I will proof a generalization.

\subsection{Properties of Sonneveld Spaces} 
Since the goal of this paper is to use Sonneveld spaces in recycling algorithms, \I first investigate their properties in more detail. 

\paragraph{Dimension and General Recursion}
The following lemma will be helpfull in a subsequent section for simpler understanding of \IDRstab.
\begin{Lemma}[General Recursion of Sonneveld Spaces]
	For all Sonneveld spaces $\cG_i$, $\cG_j$, $\forall\ i<j \in \N_0$, it holds
	\begin{align}
		\cG_{j} = p_{i,j}(\bA) \cdot \big(\cG_{i}\cap\cK^\perp_{j-i}(\bA^{H};\cP)\big)\,. \label{eqn:genRecSon}
	\end{align}
\end{Lemma}
The result is trivial. To \my best knowledge this formulation of the recursion is new. For $i=0$, e.g. shown by Gutknecht \cite{IDR-Gutkn}, i.e.
\begin{align}
	\cG_j = (\bI - \omega_1 \cdot \bA) \cdot ... \cdot (\bI - \omega_j \cdot \bA) \cdot \cK_j^\perp(\bA^H ; \cP)\,, \label{eqn:GutG}
\end{align}
\we directly see
\begin{align}
	\dim(\cG_j) \leq N - \dim\big(\cK_j(\bA^H;\cP)\big)\,. \label{eqn:dimG}
\end{align}
\gI speak of the \textit{canonical} case iff the following condition
\begin{align}
	\dim\big(\cK_j(\bA^H;\cP)\big) = \min\lbrace N,\,j \cdot \dim(\cP)\rbrace \label{eqn:DimP}
\end{align}
holds \cite[p.1038, l.13ff]{IDR-report} (in an earlier version of this, accessible in \cite{IDRweb}, it is shown that a random choice of $\cP$ satisfies \eqref{eqn:DimP} with  $\approx 100\%$ probability), also referred to as \textit{generic} \cite{IDR-report} or \textit{regular} \cite[p.8]{IDR-Gutkn} case. 


\paragraph{Related Testspaces}
When a residual is restricted into a Sonneveld space, i.e. $\br \in \cG_j$, this has the same meaning as imposing certain orthogonality properties on the residual, i.e. $\br \perp \tcC_j$ for some \textit{test space} $\tcC_j$. The relation $\cG_j = \tcC_j^\perp$ is easily seen, so Sonneveld spaces are orthogonal complements of their test spaces. From $\cG_j \subseteq \cG_i\ \forall\,j>i$ it follows $\tcC_j \supseteq \tcC_i$.

From \eqref{eqn:GutG} \we see by rewriting
\begin{align}
	\begin{split}
		\cG_j &= \cK_j^\perp\Big(\bA^H;\big(p_{0,j}(\bA^H)\big)^{-1} \cdot \cP \Big)\\
		\Leftrightarrow\ \tcC_j &= \cK_j\Big(\bA^H;\big(p_{0,j}(\bA^H)\big)^{-1} \cdot \cP \Big)
	\end{split}\,,
	\label{eqn:SonKryl}
\end{align}
cf. \cite{Simoncini}, that Sonneveld's test spaces have the structure of a Krylov subspace. However from \eqref{eqn:SonKryl} \we do not directly see $\tcC_j \supseteq \tcC_i$ because the starting space of the Krylov subspace changes with $j$ due to the stabilization polynomial. It might be useful for further study to know which directions are added onto $\tcC_{j-1}$ in specific to obtain $\tcC_j$. In the following lemma \I provide the answer.
\begin{Lemma}[Sonneveld's Test Spaces]
	With the definitions from above \we have $\tcC_0 = \lbrace \bO \rbrace$ and
	\begin{align}
		\tcC_j &= \operatornamewithlimits{span}_{k=1,\dots,j}\Big\lbrace \big(p_{0,k}(\bA^H)\big)^{-1} \cdot \cP \Big\rbrace\,, \label{eqn:tcC_add}\\
		\tcC_{j} &= \tcC_{j-1} + (\bI - \omega_j \cdot \bA^H)^{-1} \cdot \tcC_{j-1}\quad \forall \, j\in \N\,. \label{eqn:tcC_Rec}
	\end{align}
	\underline{Proof:}\\
	First \I proof \eqref{eqn:tcC_add}. Let $\bx \in \cG_j$, i.e. $\bx \perp \cK_j\Big(\bA^H ; \prod_{t=1}^j(\bI - \omega_t \cdot \bA^H)^{-1} \cdot \cP \Big)$, cf. to \eqref{eqn:SonKryl} or \cite{Sleijpen1}. \We can rewrite this to
	\begin{align*}
		\bx \perp (\bA^H)^k \cdot \prod_{t=1}^j (\bI - \omega_t \cdot \bA^H)^{-1} \cdot \cP\quad \forall\, k=0,...,j-1\,.
	\end{align*}
	\We obtain equivalent conditions on $\bx$ when \we replace the powers $(\bA^H)^k$ by other polynomials of $\bA^H$ of degree $k$. Below \I use the stabilization polynomials $p_{0,k}$.
	\begin{align*}
		\bx \perp \prod_{q=1}^k (\bI - \omega_q \cdot \bA^H) \cdot \prod_{t=1}^j (\bI - \omega_t \cdot \bA^H)^{-1} \cdot \cP\quad \forall\, k=0,...,j-1\,.
	\end{align*}
	Note that in the above expression all matrices in the product commute. Dropping factors, from this \we obtain
	\begin{align}
		\bx \perp \operatornamewithlimits{span}_{k=0,\dots,j-1} \Big\lbrace \prod_{t=1}^{j-k} (\bI - \omega_{j+1-t} \cdot \bA^H)^{-1} \cdot \cP\Big\rbrace\,. \label{eqn:Lemm2_order}
	\end{align}
	From \eqref{eqn:GutG} we see that $\cG_j$ and therefore $\tcC_j$ are invariant under different orderings of $\omega_1,...,\omega_j$. 
	Thus, turning around the ordering in \eqref{eqn:Lemm2_order} (replacing $j+1-t$ by $t$) it holds
	\begin{align}
		\operatornamewithlimits{span}_{k=0,...,j-1}\Big\lbrace \big(p_{k,j}(\bA^H)\big)^{-1} \cdot \cP \Big\rbrace = \operatornamewithlimits{span}_{k=1,...,j}\Big\lbrace \big(p_{0,k}(\bA^H)\big)^{-1} \cdot \cP \Big\rbrace\,. \label{eqn:Contained_Span}
	\end{align}
	Inserting \eqref{eqn:Contained_Span} into \eqref{eqn:Lemm2_order} yields \eqref{eqn:tcC_add}. From \eqref{eqn:tcC_add} \we see
	\begin{align}
		\tcC_j &= \tcC_{j-1} + \big(p_{0,j}(\bA^H)\big)^{-1} \cdot \cP \quad \forall\,j\in\N\,. \label{eqn:C_poly_update}
	\end{align}
	From \eqref{eqn:Contained_Span} \we see that the recursion \eqref{eqn:tcC_Rec} only adds the second term of \eqref{eqn:C_poly_update} to the next level test space because the other added directions are already contained in $\tcC_{j-1}$. q.e.d.	
\end{Lemma}
\begin{Remark}[Inverses of Stabilization Polynomials]
	For simplicity in this presentation \I assume that the relaxations $\omega_j$ $\forall j \in \N$ differ from the eigenvalues of $\bA$, thus for regular $\bA$ all the above inverses exist.
\end{Remark}

\subsection{Sonneveld Spaces in the \BiCGstab Algorithm}
\IDReins and \BiCGstab are mathematically equivalent iterative methods for the solution of single linear systems. They work on the principle of restricting a residual $\br$ of a linear system iteratively into higher level Sonneveld spaces. In the following \I review their principle.

For the Sonneveld spaces the cut-space $\cP = \operatorname{span}\lbrace \bp \rbrace$ is chosen, where $\bp \in \C^N$ is just a fixed arbitrary vector\footnote{To avoid any confusion, \I do not use the terminus \textit{shadow residual}, as a residual is usually a vector that is minimized during an iterative process.}. The relaxations are chosen such that the length of the residual is locally minimized. To restrict the residual $\br$ into successively higher Sonneveld spaces both methods use an \textit{auxiliary vector} $\bv \in \C^N$.

\begin{figure}
	\centering
	\includegraphics[width=0.8\linewidth]{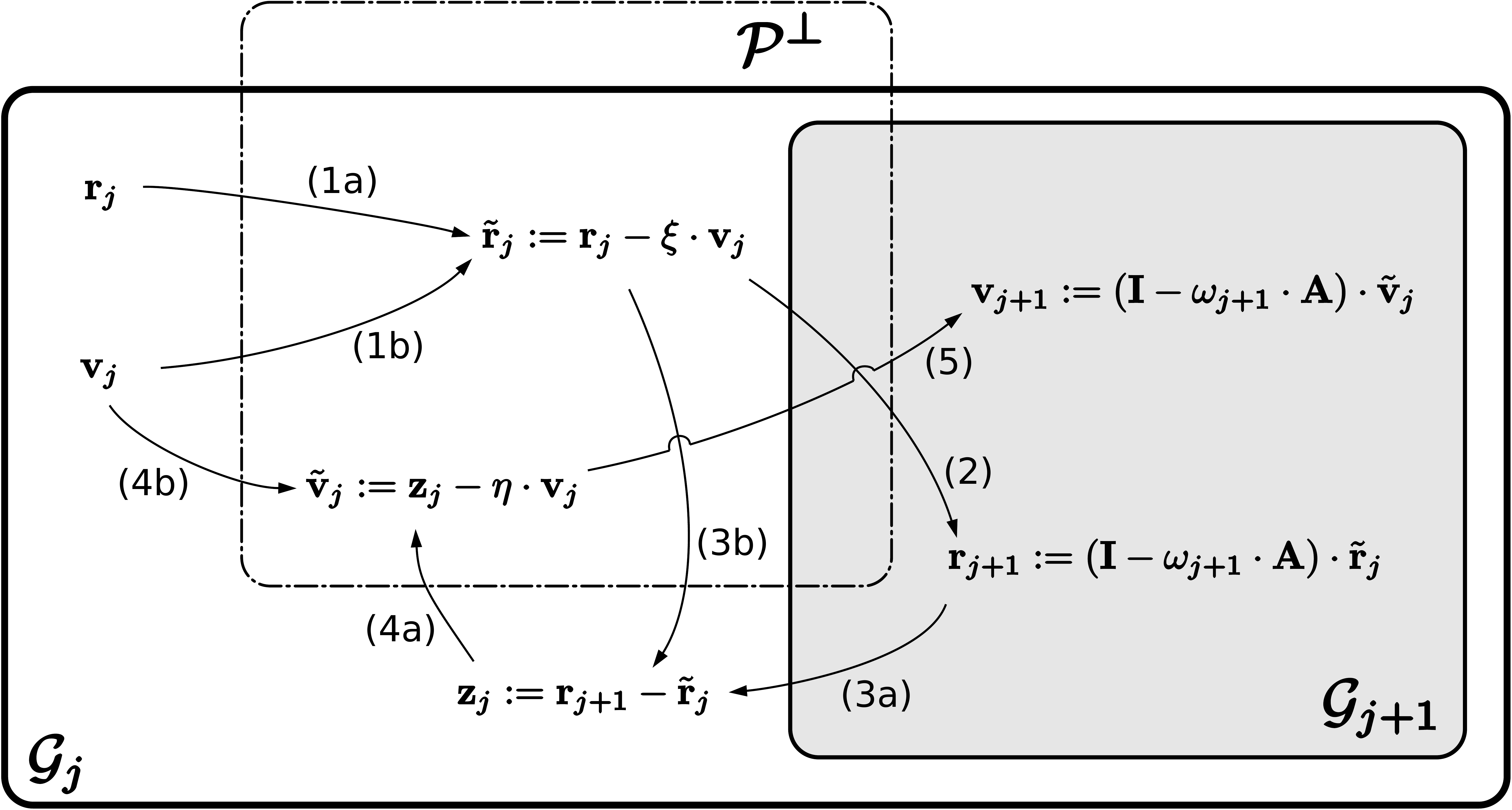}
	\caption{\textit{Computational steps of one \IDRO-cycle.}}
	\label{fig:IDRcycle}
\end{figure}

For better understanding of the iterative process \I add the iteration index $j$ as footnote to $\br,\bv$. With the vectors $\br_j,\bv_j \in \cG_j$, one iteration of both \BiCGstab and \IDReins, in the remainder called \textit{IDR cycle}, consists of the following five steps, as shown in Fig.~\ref{fig:IDRcycle}.
\begin{enumerate}
	\item Choose $\xi \in \C$, such that  $\tbr_j$ is orthogonal to $\bp$, thus $\tbr_j \in \cG_j \cap \cP^\perp$. This is sometimes called the \textit{\BiCG step}.
	\item Choose $\omega_{j+1} \in \Cno$, such that  $\br_{j+1}$ preferably has small length. Notice that $\br_{j+1} \in \cG_{j+1}$. This is sometimes referred to as \textit{\textsf{GMRES(1)} step} because a minimal residual polynomial of first order is used to minimize the length of $\br_{j+1}$.
	\item Compute $\bz_j$. Notice that $\bz_j \in \opspan\lbrace \tbr_j,\br_{j+1} \rbrace \subset \cG_{j}$.
	\item Choose $\eta \in \C$, such that  $\tbv_j$ is orthogonal to $\bp$, thus $\tbv_j \in \cG_j \cap \cP^\perp$.
	\item Compute $\bv_{j+1} \in \cG_{j+1}$.
\end{enumerate}
Using these steps, both the residual and the auxiliary vector are moved from $\cG_{j}$ to $\cG_{j+1}$, i.e. for each \IDRO-cycle the level $j$ of the Sonneveld space increases by one. In theory after sufficiently many iterations it will be $\cG_{j+1} = \lbrace \bO \rbrace$ \cite[theo. 2.1 (ii)]{IDR-report} and the residual becomes $\bO$. However, in practise the procedure often can be stopped much earlier due to sufficiently small $\|\br\|$.

\subsection{Towards \IDRstab}
There are two ways in which \BiCGstab resp. \IDReins can be improved.

\subsubsection{Improving Stability}
The first option is to generalize the {\textsf{GMRES(1)} step}. By including a {pre-calculation} stategy one can probably find better relaxations that lead to faster decreasing residual norms in the long term.

One possible approach to implement a {pre-calculation} stategy is to replace the recursion \eqref{eqn:DefSonneveld} by \eqref{eqn:genRecSon} to move from $\cG_j$ to $\cG_{j+\ell}$. In this the polynomial $p_{j,j+\ell}$ can be {customized} for residual minimization. \We review this in the following. 

\begin{figure}
	\centering
	\includegraphics[width=0.6\linewidth]{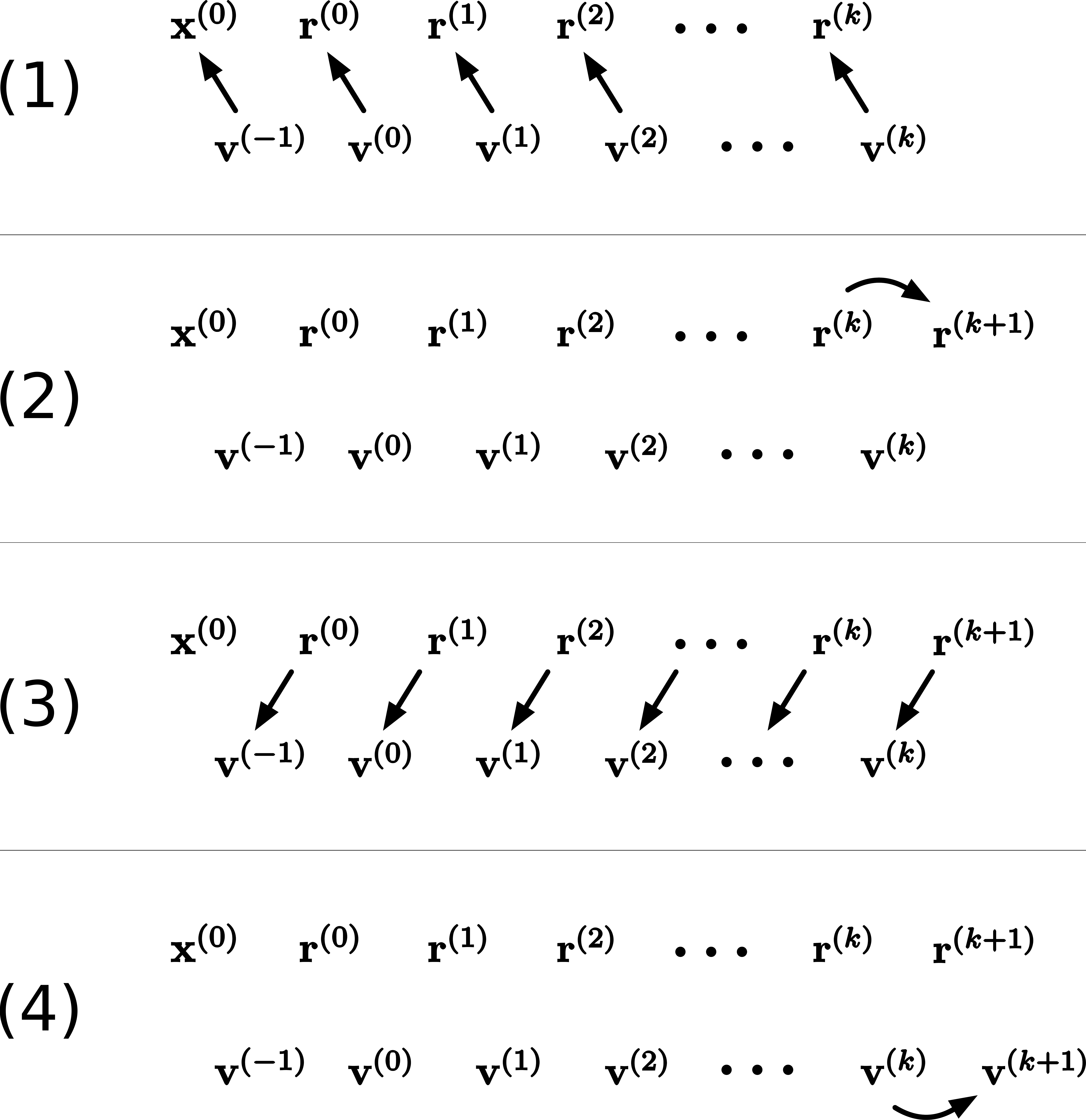}
	\caption{\textit{Computational steps of one level iteration for the pre-calculation strategy.}}
	\label{fig:LookAhead}
\end{figure}

Fig.~\ref{fig:LookAhead} shows so called \textit{level vectors} of both the residual and the auxiliary vector.
\begin{Definition}[Level Vector]
	Let $\bg \in \C^N$ be an arbitrary vector. $\bg^{(k)} := \bA^k \cdot \bg$ is called the $k^\text{th}$ level vector of $\bg$.
\end{Definition}

\We are starting from $\br\equiv \br^{(0)},\bv\equiv \bv^{(0)} \in \cG_j$ with corresponding numerical solution $\bx^{(0)}$, i.e. $\br^{(0)} = \bb - \bA \cdot \bx^{(0)}$, and pre-image $\bv^{(-1)} = \bA^{-1} \cdot \bv^{(0)}$. The following \textit{level iteration} shows how to compute level vectors $\br^{(k)},\bv^{(k)}$ of $\br,\bv$ for $k=0,...,\ell$, such that  afterwards additionally $\br,\bv \in \cK^\perp_\ell(\bA^H;\cP)$ holds. The iteration goes for $k=0,...,\ell-1$ and consists of four steps as shown in Fig.~\ref{fig:LookAhead}:
\begin{enumerate}
	\item Given $\br^{(k)},\bv^{(k)}$, orthogonalize $\br^{(k)}$ with $\bv^{(k)}$ onto $\bp$:
	\begin{align*}
		\br^{(k)} := \br^{(k)} - \xi \cdot \bv^{(k)}
	\end{align*}
	To maintain the level vector properties, use the same update for the pre-images of $\br^{(k)}$ and $-$ with sign-flip $-$ for the numerical solution $\bx^{(0)}$ of $\br^{(0)}$.
	\item Compute $\br^{(k+1)}$ from $\br^{(k)}$.
	\item Orthogonalize $\bv^{(k)}$ with $\br^{(k+1)}$ onto $\bp$, and preserve the level vector property.
	\begin{align*}
		\bv^{(i)} := \bv^{(i)} - \eta \cdot \br^{(i+1)} \quad \forall\,i=-1,...,k\,.
	\end{align*}
	\item Compute $\bv^{(k+1)}$ from $\bv^{(k)}$.
\end{enumerate}
The use of the computed data from the level iterations is as follows: \We have $\br^{(0)},\bv^{(0)} \in \cG_j \cap \cK_\ell^\perp(\bA^H;\cP)$. Choosing a stabilization polynomial
\begin{align*}
p_{j,j+\ell}(t)=\prod_{k=1}^{\ell}(1-\omega_{j+k} \cdot t) = 1 - \sum_{k=1}^\ell \tau_k \cdot t^k\,,
\end{align*}
\we can construct a new residual and auxiliary vector in the next Sonneveld space from the level vectors without any new computations of matrix-vector-products:
\begin{align*}
	\br &:= p_{j,j + \ell}(\bA) \cdot \br^{(0)}\equiv \br^{(0)} - \sum_{k=1}^\ell \tau_k \cdot \br^{(k)} \in \cG_{j+\ell}\\
	\bv &:= p_{j,j + \ell}(\bA) \cdot \bv^{(0)}\equiv \bv^{(0)} - \sum_{k=1}^\ell \tau_k \cdot \bv^{(k)} \in \cG_{j+\ell}
\end{align*}
In this \we choose the values for $\tau_1,...,\tau_\ell$ such that $\br \perp \lbrace \br^{(1)},...,\br^{(\ell)}\rbrace$, i.e. $\|\br\|_2$ is minimized. The actual relaxations $\omega_{j+1},...,\omega_{j+\ell}$ can be left unknown.
\largeparbreak
This is how higher order stabilization can be achieved. The resulting method is called \BiCGstabL \cite{BiCGstabL}. The pre-calculation strategy increases the stability and is very usefull for highly asymmetric systems. An implementation is given at a later moment.

\subsubsection{Improving Efficiency}
So far one auxiliary vector $\bv$ was used to orthogonalize the residual $\br$ onto the one-dimensional space $\cP = \opspan\lbrace\bp\rbrace$. Though the recursive definition of Sonneveld spaces also allows for higher dimensional cut-spaces $\cP$. Thus by use of $s$ auxiliary vectors $\bv_1,...,\bv_s \in \C^N$ and $s$ basis vectors $\bp_1,...,\bp_{s} \in \C^N$ for $\cP$ instead of one, the iteration from $\br_j,\bv_{j,1},...,\bv_{j,s} \in \cG_{j}$ to $\cG_{j+1}$ now works as follows (cf. Fig.~\ref{fig:IDRcycle}):
\begin{Aufz}
	\item Orthogonalize $\br_j$ with $\bv_{j,1},...,\bv_{j,s}$ onto $\cP$ to obtain $\tbr_j$.
	\item Compute $\bz_j = \bA \cdot \tbr_j$.
	\item Choose $\omega_j$ such that  the next residual is minimized and then compute it:
	\begin{align*}
		\br_{j+1} := \tbr_j - \omega_j \cdot \bz_j\quad //\, \br_{j+1} \in \cG_{j+1}\,.
	\end{align*}
	\item Then compute for all auxiliary vectors $k=1,...,s$ in order:
	\begin{Aufz}
		\item Orthogonalize $\bv_{j,k}$ with the other auxiliary vectors and $\bz_j$ onto $\cP$ to obtain $\tbv_{j,k}$:
		\begin{align*}
			\tbv_{j,k} := \bz_j - \sum_{q=1}^{k-1} \eta_{k,q} \cdot \bv_{j+1,q} - \sum_{q=k}^{s} \eta_{k,q} \cdot \bv_{j,q}\quad // \, \tbv_{j,k} \in \cG_{j-1} \cap \cP^\perp\,.
		\end{align*}
		\item Move $\tbv_{j,k}$ from $\cG_{j}$ into $\cG_{j+1}$:
		\begin{align*}
			\bv_{j+1,k} := \tbv_{j,k} - \omega_j \cdot \bA \cdot \tbv_{j,k}\quad // \, \bv_{j+1,k} := (\bI - \omega_j \cdot \bA) \cdot \tbv_{j,k} \in \cG_{j+1}\,.
		\end{align*}
	\end{Aufz}
\end{Aufz}
So in principle the computational steps of the \IDRO-cycle stay the same as in Fig.~\ref{fig:IDRcycle}; one only uses the fact that one vector can be orthogonalized with $s$ others against $s$ arbitrary space directions.

The use of the larger dimension $s$ of $\cP$ is that $-$ only using $s+1$ matrix-vector-products (\MV) with $\bA$ $-$ the dimension of the Sonneveld space, in which the residual lies, decreases by $s$. So the average ratio between computed matrix-vector-products and dimensions of the test space (remember $\tcC_j$) is $\lambda = s/(s+1)$ and approaches the optimal\footnote{for the case of solving one single system. \We see later that for sequences of systems there is no limit for $\lambda$.} value $1$ for moderate values of $s$, e.g. $s=4$. \gI call $\lambda$ \textit{efficiency ratio} because it gives a ratio of reduced dimensions over computational effort. In practical numerical applications it is observed that $\lambda$ correlates with superlinearity, i.e. higher values of $s$ increase the superlinear convergence. Moreover a higher $s$ can sometimes improve the stability of the iterative method.
\largeparbreak
Methods that use the above concept are \IDR \cite{IDR-report} and \textsf{ML($s$)BiCGStab} \cite{MLBiCGstab}. $s$ denotes $\dim(\cP)$ and the number of auxiliary vectors. The vectors $\bp_1,...,\bp_s$ can be chosen arbitrarily, e.g. randomly.

Methods that use both multiple auxiliary vectors (i.e. $s > 1$) and the pre-calculation strategy (i.e. $\ell > 1$) are \IDRstab \cite{IDRstab-Paper} and \textsf{GBi-CGStab($s,\ell$)} \cite{GBiCGstab}. 

\paragraph{Implementation}
Algorithm \ref{algo:IDRstab} gives an implementation of \IDRstab with comments and emphasis on simplicity. Robust and preconditioned variants for practical use are compared and can be found in \cite{IDRstabS2} and \cite{Aihara2}.

\begin{algorithm}
	\caption{\IDRstab}
	\label{algo:IDRstab}
	\begin{algorithmic}[1]
		\Procedure{IDRstab}{$\bA,\bb,\bx,s,\ell,\cP,\bU,\bV$}
		\State \textit{//prior requirements: $\bV = \bA \cdot \bU \in \C^{N \times s}$, $\dim(\cP) = s$}
		\State $\br := \bb - \bA \cdot \bx$\,,\quad$j:=0$
		\While{$\|\br\|>\tol$}
		\State $\br^{(0)} := \br$, $\bx^{(0)} := \bx$
		\State $\bv_q^{(-1)} := \bu_q$, $\bv_q^{(0)} := \bv_q$\quad$\forall\,q=1,...,s$\quad \textit{//$\bU = [\bu_1,...,\bu_s],\,\bV = [\bv_1,...,\bv_s]$}
		\State \textit{// \IDRO-steps / level iteration}
		\For{$k=0,1,...,\ell-1$}
		\State Choose $\boldsymbol{\gamma} \in \C^s$, such that  $\br^{(k)} - \bV^{(k)} \cdot \boldsymbol{\gamma} \perp \cP$.
		\State \textit{// notation $\bV_{d:f}^{(i)} \equiv [\bv_d^{(i)},...,\bv_f^{(i)}]$}
		\State $\br^{(i)} := \br^{(i)} - \bV^{(i)} \cdot \boldsymbol{\gamma}$\quad, $\forall\,i=0,1,...,k$\,.
		\State Recapture solution to $\br^{(0)}$ : $\bx^{(0)} := \bx^{(0)} + \bV^{(-1)} \cdot \boldsymbol{\gamma}$\,.
		\State $\br^{(k+1)} := \bA \cdot \br^{(k)}$
		\For{$q=1,2,...,s$}
		\State Choose $\boldsymbol{\eta}_q \in \C^s$, such that  $\br^{(k+1)} - [\bV_{1:q-1}^{(k+1)},\bV_{q:s}^{(k)}] \cdot \boldsymbol{\eta}_q \perp \cP$.
		\State $\bv_q^{(i)} := \br^{(i+1)} - [\bV_{1:q-1}^{(i+1)},\bV_{q:s}^{(i)}] \cdot \boldsymbol{\eta}_q$\quad$\forall\, i=-1,0,...,k$\,.
		\State $\bv_{q}^{(k+1)} := \bA \cdot \bv_q^{(k)}$
		\EndFor
		\EndFor
		\textit{// $\br^{(0)},\bv_1^{(0)},...,\bv_s^{(0)} \in \cG_j \cap \cK^\perp_\ell(\bA^H;\cP)$}
		\State \textit{// degree $\ell$ residual minimization}
		\State $\bZ := [\br^{(1)},...,\br^{(\ell)}] \in \C^{N \times \ell}$
		\State $\boldsymbol{\tau} \equiv (\tau_1,...,\tau_\ell)^T := \bZ^\dagger \cdot \br^{(0)}$
		\State $\br := \br^{(0)} - \sum_{i=1}^{\ell} \tau_i \cdot \br^{(i)}$\quad \textit{// $\equiv \br^{(0)} - \bZ \cdot \boldsymbol{\tau}$} \label{algo1:Z23}
		\State Recapture solution: $\bx := \bx^{(0)} + \sum_{k=i}^{\ell} \tau_i \cdot \br^{(i-1)}$ \label{algo1:Z24}
		\State $\bv_q := \bv_q^{(0)} - \sum_{i=1}^{\ell} \tau_i \cdot \bv_q^{(i)}$\quad$\forall\, q=1,...,s$\,.\label{algo1:Z25}
		\State $\bu_q := \bv_q^{(-1)} - \sum_{i=1}^{\ell} \tau_i \cdot \bv_q^{(i-1)}$\quad$\forall\, q=1,...,s$\,.\label{algo1:Z26}
		\State $j := j + \ell$
		\EndWhile
		\State $J := j$
		\State \Return $J$,$\cP,\bU,\bV,\bx,\br$
		\EndProcedure
	\end{algorithmic}
\end{algorithm}

\section{Recycling in \IDRO Methods}

Early \IDRO-based Krylov subspace recycling methods\footnote{i.e. Miltenberger's method and \SRIDR} are based on \IDR without the stabilization approach, thus $\ell = 1$, and the relaxation parameters $\omega_j\ \forall j \in \N$ are explicitly known.

\subsection{Miltenberger's Approach}
In \cite[chap. 4.3]{IDR-Milten} Miltenberger proposes an early \IDRO-method with a recycling technique. For the solution of sequences like \eqref{eqn:AxbSequence} he suggests to push the result matrices $\bU,\bV$ of algorithm \ref{algo:IDRstab} from the solution process of a \rhs $\bb\hia$ as input arguments for the solution process of a subsequent system with \rhs $\bb^{(\iota+1)}$. To give an example:

A first system $\bA \cdot \bx\ha = \bb\ha$ is solved with e.g. $\bU = [\be_1,...,\be_s]$, $\bV = \bA_{1:s}$ and \IDR. From the solution process besides to $\br$ and $\bx$ two overwritten matrices $\bU,\bV$ are returned. \textit{These} two matrices are now used as input arguments for \IDR for the solution of the subsequent \rhs $\bb\hb$.

This is a reasonable first recycling approach. However, at first glance one would not expect that this method recycles more than the $s$-dimensional search space that is spanned by $\bU$. However, \we will see later that it actually does.

\subsection{The \SRIDR Method}
For the outputs $\bU,\bV$, that Miltenberger reuses, \we obviously have $\opImage(\bV) \subseteq \cG_J$, where $J \in \N$ is the number of \IDRO-cycles that were performed in the former solution process. The specific $\cG_J$ is the last Sonneveld space in which the residual of the solution process with outputs $\bU,\bV$ was restricted. Continuing Miltenberger's approach, \I reused in \SRIDR the returned basis matrices $\bU,\bV$, to restrict an arbitrary residual of a subsequent \rhs inexpensively into that specific $\cG_J$. In the following the idea for this is reviewed more in detail.

\begin{figure}
	\centering
	\includegraphics[width=1\linewidth]{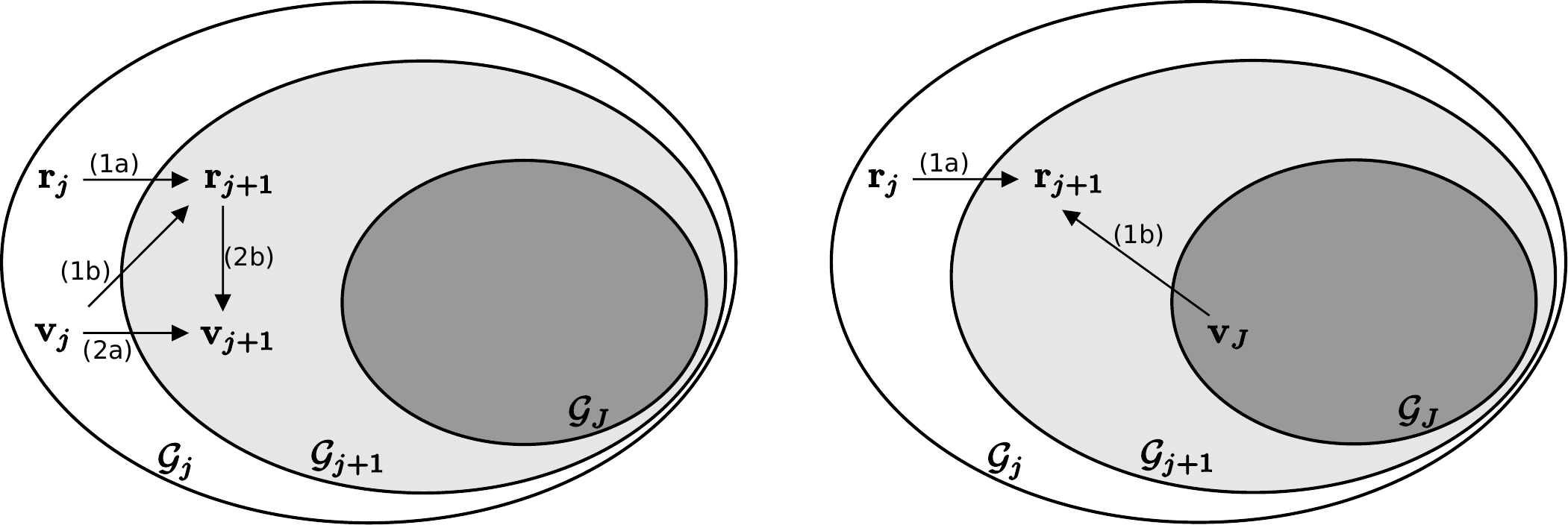}
	\caption{\textit{IDR cycle (left) compared to SRIDR cycle (right).}}
	\label{fig:SRIDRcycle}
\end{figure}

\We first consider the case $s=\ell=1$. In Fig.~\ref{fig:SRIDRcycle} left \I redraw in a compacter form the usual \IDRO-cycle as sketched in Fig.~\ref{fig:IDRcycle}: Given $\br_j,\bv_j \in \cG_j$, \we first use $\br_j,\bv_j$ to construct $\br_{j+1} \in \cG_{j+1}$ by effort of one \MV with $\bA$. Afterwards $\br_{j+1},\bv_j$ are used to compute $\bv_{j+1} \in \cG_{j+1}$ by effort of one \MV with $\bA$, too. These two steps are shown in the figure by black arrows.

In the right part of Fig.~\ref{fig:SRIDRcycle} \we have a \textit{modified} \IDRO-cycle, where instead of $\bv_j \in \cG_j$ \we use a vector $\bv_J \in \cG_J \subseteq \cG_j$. In the construction of $\br_{j+1}$ \we can replace $\bv_j$ by $\bv_J$ and still obtain $\br_{j+1} \in \cG_{j+1}$. As an advantage, the second \MV with $\bA$ for the auxiliary vector $\bv_J$ can be skipped, whenever already $\cG_J \subseteq \cG_j$, i.e. $J \geq j$. This increases the efficiency ratio from $\lambda = 1/2$ for an \IDRO-cycle to $\lambda = 1$ for a modified \IDRO-cycle.

For $s=1$, the matrix $\bV$ has only one column, which is indeed the auxiliary vector $\bv_J$. The concept of modified \IDRO-cycles can be generalized for higher values of $s$, i.e. multiple auxiliary vectors in $\cG_J$ are used in each cycle to orthogonalize $\br_j$ against multiple test vectors $\bp_1,...,\bp_s$ of a $s$-dimensional cut-space $\cP$. By using such modified \IDRO-cycles, the residual of a subsequent \rhs can be restricted to the Sonneveld space $\cG_J$ from a former solution process with a computational effort that is dominated by $J$ \MV-s with $\bA$. The efficiency ratio for the modified \IDRO-cycles is then $\lambda = s / 1 = s$ in the canonical case, i.e. validity of \eqref{eqn:DimP}.
\largeparbreak
The idea of the method \SRIDR is as follows. By reusing matrices $\bU,\bV$ with $\opImage(\bV) \subseteq \cG_J$ from a former solution process, a residual of a subsequent system is inexpensively restricted from $\cG_0$ into that former Sonneveld space $\cG_J$. Optionally, afterwards usual \IDRO-cycles are performed to restrict the residual from $\cG_J$ into even higher level Sonneveld spaces to further reduce $\|\br\|$. \SRIDR stands for \textit{short recycling} for \IDR. Short recycling is adapted to \textit{short recurrence}, as like short recurrence methods \SRIDR needs only a limited storage amount to recycle subspace information of a Sonneveld space of (theoretical) arbitrary level.

\subsection{The Issue with \SRIDR}
Despite the high efficiency ratio of \SRIDR, there is still one big drawback: For the recycling of Sonneveld spaces from the former solution process, the original relaxation parameters $\omega_1,...,\omega_J$ from the former solution process must be reused, too. Otherwise the Sonneveld spaces for our new residual would differ and the recycled auxiliary vectors would loose their subrange property. I.e. in Fig.~\ref{fig:SRIDRcycle} \we could no longer ensure $\cG_J \subseteq \cG_{j+1}$ because a change of $\omega_{j+1}$ leads to a different $\cG_{j+1}$. As the relaxations must be kept fixed during the recycling iterations, \we lack in a residual-minimizing property during this phase.

In fact for a naive implementation of \SRIDR \we see (e.g. \cite[Fig. 6]{Report15}) that during the first $J$ iterations, in which the residual is moved from $\cG_0$ to $\cG_J$ with modified \IDRO-cycles, the residual norm increases dramatically. This is probably due to the lack of a residual-minimizing property. However, \we also see from the same figure that the restiction of a residual into a high level Sonneveld space leads to a faster convergence \textit{after} the recycling iterations.
\largeparbreak
In summary, for \SRIDR there is a trade-off: On the one hand arbitrary large test spaces can be recycled by only storing a few auxiliary vectors\footnote{at least in theory where no round-off occurs}, which may improve the convergence. On the other hand, the recycling iterations lack a residual-minimizing property which may cause a blow up in the residual norm.
\largeparbreak
It would be desirable to have a method that has all the advantages of \SRIDR, i.e. achieves the same efficiency as \SRIDR, and moreover allows for a completely free choice of the relaxations. \gI present exactly such a method in the next section.

\section{A Generalization of Sonneveld Spaces suited for Reuse of Subspace Information}
In this section \we will see that Miltenberger's method has comparable termination properties to \SRIDR but allows a free choice of the relaxations. This motivates a generalization of Miltenberger's method, which is indeed \Mstab.

\subsection{Generalization of Sonneveld Spaces}
\gI introduce a generalization of Sonneveld spaces. The benefit of this generalization is that one can add arbitrary directions to a sequence of subspaces without destroying their recursive relations.
\begin{Definition}[$\cM$-space]\textcolor{white}{.}\\
	Let $\bA \in \C^{N \times N}$ be a matrix, $\lbrace \omega_j \rbrace_{j \in \N} \subset \Cno$ a sequence, $\lbrace\cQ_j\rbrace_{j \in \N} \subseteq \C^N\,,\ \lbrace\cP_j\rbrace_{j \in \N} \subseteq \C^N$ sequences of subspaces of $\C^N$ with $\cQ_{j+1} \subseteq \cQ_{j}\,,\ \cP_{j+1} \supseteq \cP_j$ $\forall j \in \N$, $\cM_0 = \C^N$ and recursively
	\begin{align*}
		\cM_j = (\bI - \omega_j \cdot \bA) \cdot (\cM_{j-1} \cap \cP_j^\perp )+\cQ_j \quad\ \forall\ j \in \N\,.
	\end{align*}
	The spaces $\cM_j$ are called \emph{$\cM$-spaces}, the $\cQ_j$ are \emph{add-spaces} and the $\cP_j$ are called \emph{cut-spaces}. $j$ is called \emph{level}.
\end{Definition}
\largeparbreak
\begin{Theorem}[Nestedness of  $\cM$-spaces]\label{theo:Mspace}\textcolor{white}{.}\\
	Let $\lbrace\cQ\ha_j\rbrace_{j \in \N},\lbrace\cQ\hb_j\rbrace_{j \in \N}$ and $\lbrace\cP\ha_j\rbrace_{j \in \N},\lbrace\cP\hb_j\rbrace_{j \in \N}$ be two sequences of add- and cut-spaces respectively, with $\cQ\ha_j \subseteq \cQ\hb_j$, $\cP\ha_j \supseteq \cP\hb_j$ $\forall j \in \N$. Let $\cM\hia_j$, $j \in \N_0$ be the $\cM$-spaces for $\iota=1,2$ respectively.
	
	Then the following holds:
	\begin{align*}
		\cM\hia_{j+1} &\subseteq \cM\hia_j\quad \forall j \in \N_0\,,\ \iota \in \lbrace 1,2\rbrace\\
		\cM\ha_{j} &\subseteq \cM\hb_j\quad \forall j \in \N_0\,.
	\end{align*}
	\underline{Proof:}\\
	\gI show the first result for $\iota = 1$ and drop the super index.\\
	By Induction:
	\begin{Aufz}
		\item Basis: Obviously $\cM_1 \subseteq \cM_0$.
		\item Hypothesis: $\cM_j \subseteq \cM_{j-1}$ holds for some $j \in \N$.
		\item Induction step: $\cM_j \subseteq \cM_{j-1} \, \Rightarrow \, \cM_{j+1} \subseteq \cM_{j}$ is to show.\\
		Choose an arbitrary $\bx \in \cM_{j+1}$. Then
		\begin{align*}
			\exists\, \by \in \cM_{j} \cap \cP_{j+1}^\perp \ \wedge \ \bq \in \cQ_{j+1} \,:\, \bx = (\bI - \omega_{j+1} \cdot \bA) \cdot \by + \bq\,.
		\end{align*}
		Due to $(\cM_{j} \cap \cP_{j+1}^\perp) \subseteq (\cM_{j-1} \cap \cP_{j}^\perp)$ by induction hypothesis and nestedness of the cut-spaces, \we have
		\begin{align*}
			\tbx := (\bI - \omega_j \cdot \bA) \cdot \by \in \cM_j\,.
		\end{align*}
		With $\cQ_{j+1} \subseteq \cQ_j \subseteq \cM_j$ \we have $\by,\tbx,\bq \in \cM_j$. As $\bx \in \operatorname{span}\lbrace \by,\bq,\tbx\rbrace \subset \cM_j$ for arbitrary $\bx \in \cM_{j+1}$, it is shown that $\cM_{j+1} \subseteq \cM_j$.
	\end{Aufz}
	The second result can also be shown by induction with the same idea.\\
	q.e.d.
\end{Theorem}
\begin{Remark}[Properties of $\cM_j$]\label{rem:Mspaces}
	\textcolor{white}{.}
	\begin{enumerate}
		\item For $\cQ_1 = \lbrace \bO\rbrace$, $\cP_j = \cP$ $\forall j \in \N$ one obtains Sonneveld spaces.
		\item One can describe the recursion for $\cM_j$ by cutting (with $\cP_j^\perp$), shearing (by multiplication with shifted $\bA$) and adding (by $\cQ_j$). The order in this is \enquote{cut-shear-add}. The nestedness property of the $\cM$-spaces is conserved if this order is changed in all recursions in the same way, e.g. to \enquote{add-cut-shear}.
		\item One can expect that in the canonical case it holds
		\begin{align*}
			\dim(\cM_j) \leq \max\big\lbrace 0\,,\,\dim(\cM_{j-1}) - \dim(\cP_j)\big\rbrace + \dim(\cQ_j)\,.
		\end{align*}
		Numerical experiments indicate that this is sharp. Choosing $\cQ_j = \lbrace \bO \rbrace$ for all $j$ greater than a certain number yields $\cM_j = \lbrace \bO \rbrace$ for a sufficiently high level.
		\item In the proof the relaxations drop out by using $\by \in \cM_j$. Thus the intuition from Fig.~\ref{fig:DialPlate} is precisely the mathematical reason for the nestedness of the $\cM$-spaces and Sonneveld spaces.
	\end{enumerate}
\end{Remark}
Roughly speaking, the advantage of the $\cM$-spaces is that they have the same nestedness properties as Sonneveld spaces, but one can add arbitrary directions to them with only mild growth of their dimension. One can think of many applications where these spaces are useful, e.g. in alternative stabilization approaches.

\subsection{Application for Recycling Methods: The \Mstab Method}
As we consider sequences of linear systems, \I propose a way to use $\cM$-spaces for the numerical solution of sequences of linear systems:

Consider the case \we first solve a system $\bA \cdot \bx\ha = \bb\ha$ with a usual \IDRO-method for some $\cP = \opspan\lbrace\bp_1,...,\bp_s\rbrace$ and relaxations $\omega_1,...,\omega_J$ for some $J \in \N$. The \IDRO-method constructs the sequence $\cG_0,...,\cG_J$ of Sonneveld spaces. As a byproduct of the solution process \we obtain auxialiary vectors $\bv_1,...,\bv_s \in \cG_J$.

Next \we want to solve a subsequent system $\bA \cdot \bx\hb = \bb\hb$ with initial residual $\br_0 := \bb\hb - \bA \cdot \bx_0$. Instead of \SRIDR, where \we manipulated $\br_0$, to move it from $\cG_0$ to $\cG_J$, \we can now instead manipulate the Sonneveld space $\cG_J$ itself.
\largeparbreak
The overall strategy is to widen $\cG_J$ to a slightly larger $\cM$-space $\cM_J$, such that  it contains $\br_0$. Then by use of \IDRO-cycles on the $\cM$-spaces, this residual can be improved iteratively by shrinking $\cM_J$.

\paragraph{Construction} As a strategy, \I first construct add- and cut-spaces with superindex $1$, such that  the $\cM$-spaces of these are identical to the Sonneveld spaces. Remark \ref{rem:Mspaces} point 1 tells how this can be done. Then \I construct add- and cut-spaces with superindex $2$, such that  the requirements of theorem \ref{theo:Mspace} are satisfied, and such that  $\br_0 \in \cM\hb_J$ holds.

By choosing $\cQ\ha_j = \lbrace \bO \rbrace$ and $\cP\ha_j = \cP$ $\forall j \in \N$, we have $\cM\ha_j = \cG_j$ $\forall j\in \N_0$ for the $\cM$-spaces $\cM\ha_j$ such as defined in Theorem \ref{theo:Mspace}.

By choosing
\begin{align*}
	\cQ\hb_j = \begin{cases}
		\opspan\lbrace \br_0 \rbrace & \text{if }j\leq J\\
		\lbrace \bO \rbrace & \text{else}
	\end{cases}\quad \forall j\in\N
\end{align*}
and $\cP\hb_j = \cP$ $\forall j\in\N$, the $\cM$-spaces $\cM\hb_j$ fullfill (cf. Theorem 1)
\begin{align*}
	\cG_j = \cM\ha_j \subseteq \cM\hb_j\quad\forall j\in\N_0\,.
\end{align*}
Obviously, by our choice of the add-spaces $\cQ\hb_j$, $j=1,...,J$ \we also ensured $\br_0 \in \cM\hb_J$. In the remainder \I drop the superindices and speak of $\cG_j$ (for $\iota=1$) and $\cM_j$ (for $\iota=2$) for resp. level $j$. \We notice
\begin{align*}
	\dim(\cM_J) \leq \dim(\cG_J) + J\,,
\end{align*}
as in each recursion for $\cM_j$, $j \leq J$, only a one-dimensional space $\cQ_j$ was added.

\begin{figure}
	\centering
	\includegraphics[width=1\linewidth]{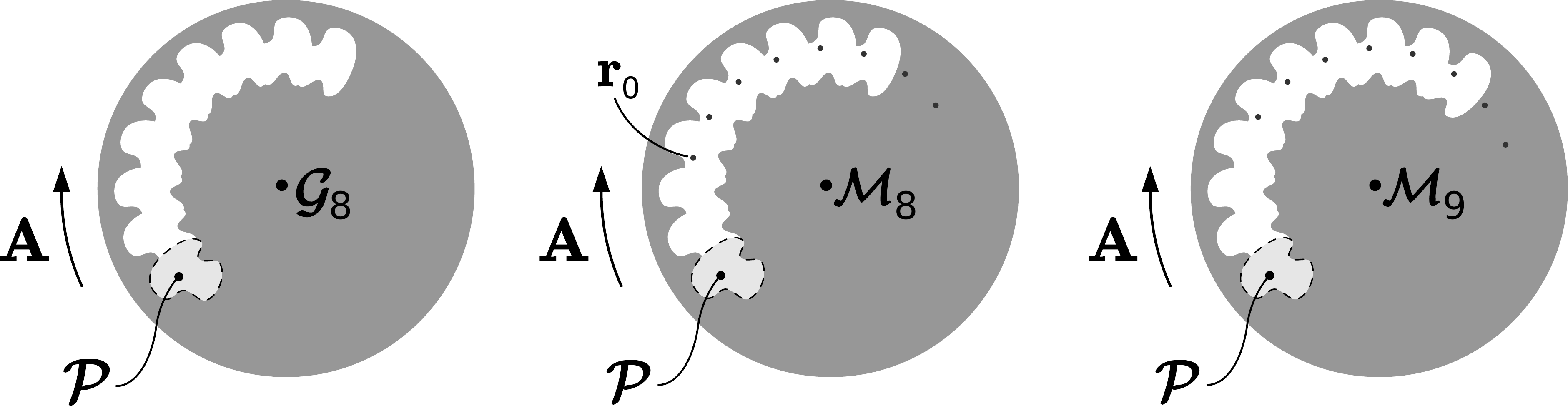}
	\caption{Geometric Intuition of $\cM$-spaces (center \& right) compared to Sonneveld spaces (left).}
	\label{fig:M_space}
\end{figure}

\paragraph{Intuition} To have an imagination of $\cM$-spaces compared to Sonneveld spaces, Fig.~\ref{fig:M_space} gives an illustration for the Sonneveld space (left) and the $\cM$-space (centre) that is obtained for $j=8$, $\cQ_1=...=\cQ_8 = \opspan\lbrace\br_0\rbrace$ for arbitrary fictional $\br_0$ and $\cP$ (this time not a sphere). Comparing $\cM_8$ to $\cG_8$, $\cM_8$ has additional dimensions that originate from level vectors of $\br_0$. In the figure these dimensions are visualized by darker points, cognoscible by their equidistant positioning on a circle line. From the illustation \we see the following: If \we perform \IDRO-cycles on $\cM_8$ (cut out $\cP$, then rotate in the direction of $\bA$), then the subsequent $\cM$-spaces are contained in $\cM_8$. Thus the $\cM$-spaces have the same nestedness properties as Sonneveld spaces. To illustrate this, Fig.~\ref{fig:M_space} shows (right) the space $\cM_9$ that \we obtain for $\cQ_9 = \lbrace \bO\rbrace$. In the picture we can see that by the rotation each level vector moves to the position of its successor, thus $\cM_9 \subseteq \cM_8$.

In Fig.~\ref{fig:M_space} $\br_0$ lies in a domain that was already cutted from $\cG_8$. This is the general case. If instead $\br_0$ would lay in $\cG_8$, then \we could have chosen $\cQ_1=...=\cQ_8 = \opspan\lbrace\bO\rbrace$. As in an algorithm \we will explicitly compute neither $\cM$-spaces nor add-spaces, \we can always assume to have chosen the optimal (i.e. smallest sufficient) add-spaces.

\paragraph{The Algorithm} Finally \we need a numerical method that iteratively moves the residual $\br \in \cM_J$ from $\cM_J$ into higher level $\cM$-spaces. To achieve this any usual \IDRO-method works\footnote{All one has to do is passing $\bU,\bV$ and $\cP$ from the output of the former solution process as inputs for the current solution process.}. This is due to the reason that, as \we have chosen $\cQ_j = \lbrace \bO \rbrace$ and $\cP_j = \cP$ for all $j>J$, the recursion for subsequent $\cM$-spaces is identical to that of Sonneveld spaces, i.e. cut out $\cP$ and rotate in direction $\bA$. For later reference such \IDRO-methods are called $\cM$-methods. As an illustration, the iterative procedure from Fig.~\ref{fig:IDRcycle} can be used without modification to iterate, e.g., in Fig.~\ref{fig:M_space} vectors from $\cM_8$ to $\cM_9$.
\largeparbreak
With the background from above \we see that Miltenberger's method, which uses \IDR, was the first $\cM$-method. In this paper \I propose to use \IDRstab instead of \IDR for the iterations. As an implementation \I refer to Algo. \ref{algo:IDRstab}, where for subsequent systems the outputs $\cP,\bU,\bV$ of a former solution process must be reused as input arguments. For later reference this method is called $\Mstab$.

\subsection{Properties of \Mstab compared to \SRIDR}
To compare both methods \I define two numbers that are counted during a computational solution process of \SRIDR and \Mstab, respectively. $\hMV$ is the number of computed matrix-vector-products with $\bA$, an important cost measure for computational effort. $\hRD$ is the number of dimensions of the test space (remember $\tcC$, defined by the $\cM$- resp. Sonneveld space) against which the residual is orthogonalized. Here $\hRD$ stands for \textit{reduced dimensions}\footnote{\gI chose a different name from \emph{dimension reduction} to indicate that this is a general property of all Krylov subspace methods.}. In the following \I will always assume the canonical case, i.e. that \eqref{eqn:DimP} holds. Both methods reuse data from a Sonneveld space $\cG_J$ of level $J$.

\gI recall that for $j\leq J$ \SRIDR reduces $s$ dimensions from the residual for each computed matrix-vector-product. For $j > J$ \SRIDR applies usual \IDRO-cycles by which in average $s/(s+1)$ dimensions are reduced per matrix-vector-product.

In \Mstab instead the fact is used that for $\cQ_{j\leq J} = \opspan\lbrace \br_0 \rbrace$ an arbitrary residual $\br_0 \in \C^N$ already lies in $\cM_J \supseteq \cG_J$, which has at most $J$ more dimensions than $\cG_J$. Thus without any computations, \Mstab already starts at $\hRD = J \cdot (s-1)$. \Mstab uses \IDRO-cycles with pre-calculations, thus reduces $s$ additional dimensions from the residual within each \IDRO-cycle.

\begin{table}[H]
	\centering
	\begin{tabular}{||c||l | l | l ||}
		\hline
		\hline
		\textbf{Method} & \textbf\hMV & \textbf\hRD & \textbf{\#Columns} \\
		\hline
		\SRIDR
		& $j \cdot (s+1) - \min\lbrace j,\,J\rbrace \cdot s$ & $j \cdot s\phantom{{} + J \cdot (s-1){}}$ & $1 + (s+1) \cdot 2$ \\
		\Mstab & $j \cdot (s+1)\phantom{{} - \min\lbrace j,\,J\rbrace \cdot s}$ & $j \cdot s + J \cdot (s-1)$ & $(s+1) \cdot (\ell+2)$ \\ 
		\GMRES & $j \cdot s$ & $j \cdot s$ & $1+s \cdot j$ \\ 
		\hline
		\hline
	\end{tabular}\\
	\textcolor{white}{\footnotesize.}
	\caption{{Matrix-vector-products and reduced dimensions for \SRIDR and \Mstab after $j$ \IDRO-cycles in the canonical case. For comparison \I added \GMRES, where $j \cdot s$ expresses the number of iterations. \#Columns is the number of stored column vectors, which is unlimited only for \GMRES.}}
	\label{tab:MVsRDs}
\end{table}
\We see from Tab. \ref{tab:MVsRDs} that \SRIDR and \Mstab have roughly the same efficiency for not to small $s$ and $j$. However, an important practical advantage of \Mstab over \SRIDR is the freedom of choosing the relaxation parameters arbitrarily. Besides \we see that both methods have a higher efficiency than \GMRES, thus they probably terminate faster (as \we will actually see for \Mstab in the numerical experiments and have seen for \SRIDR in \cite[Fig. 6]{Report15}).

\section{Numerical Experiments with \Mstab}
This section is organized as follows. First of all \I solve a test case that gives evidence to the numerical efficiency of \Mstab. Afterwards \I investigate the finite termination behaviour of \Mstab in more detail.

\paragraph{Preliminaries}
In the following examples where \I test \IDRstab and \Mstab, \I use for both methods the implementation from Algo. \ref{algo:IDRstab}. Therefore in convergence graphs the residuals of $\bx$ from line \ref{algo1:Z24} are plotted as dots and connected with lines.

When using \IDRstab, the input matrices $\bU,\bV$ are obtained by an Arnoldi scheme of level $s$. The cost of this in \hMV is accounted for by a shift of the convergence graph by $s$ positions to the right.

When using \Mstab, the input matrices $\bU,\bV$ instead are \emph{fetched} during a run of \IDRstab. To be more precise: For each dot in a convergence graph of \IDRstab a new solution $\bx$ (cf. Algo. \ref{algo:IDRstab} l. \ref{algo1:Z24}) and new matrices $\bU,\bV$ (cf. l. \ref{algo1:Z25}-\ref{algo1:Z26}) are computed. By choosing one of these dots as \emph{fetching point}, the according matrices $\bU,\bV$ are \enquote{fetched}, i.e. written out, and used as inputs for \Mstab.

\subsection{Numerical Efficiency}\label{sec:5.1}
Here a test case with a nonsymmetric sequence of linear systems is shown, where \Mstab has a superior numerical efficiency over the common Krylov subspace methods \GMRES, \BiCG, \BiCGstab and \IDRstab.

This test problem is accessible in different sizes in \cite{IDRweb}. The system results from a finite element discretization of an ocean model \cite{Ocean}. \We study the largest available test case, with $N=42249$ with $\bA \in \R^{N \times N}$, $\cond_1(\bA) \approx 6.17 \cdot 10^7$ with twelve \rhs-es $\bb\ha,\bb\hb,...,\bb^{(12)}$, resulting from month-dependent wind fields. \gI use a splitted preconditioning
\begin{align*}
	\bL^{-1} \cdot \bA \cdot \bR^{-1} \cdot \bx\hia = \bL^{-1} \cdot \bb\hia\quad \iota = 1,...,12\,,
\end{align*}
with $\bL,\bR$ resulting from an incomplete LU-factorization with zero fill-in. Each respective system is solved for $\bx\hia$. The solution to the original system is then $\bR^{-1} \cdot \bx\hia$, respectively.

\begin{figure}
\centering
\includegraphics[width=1\linewidth]{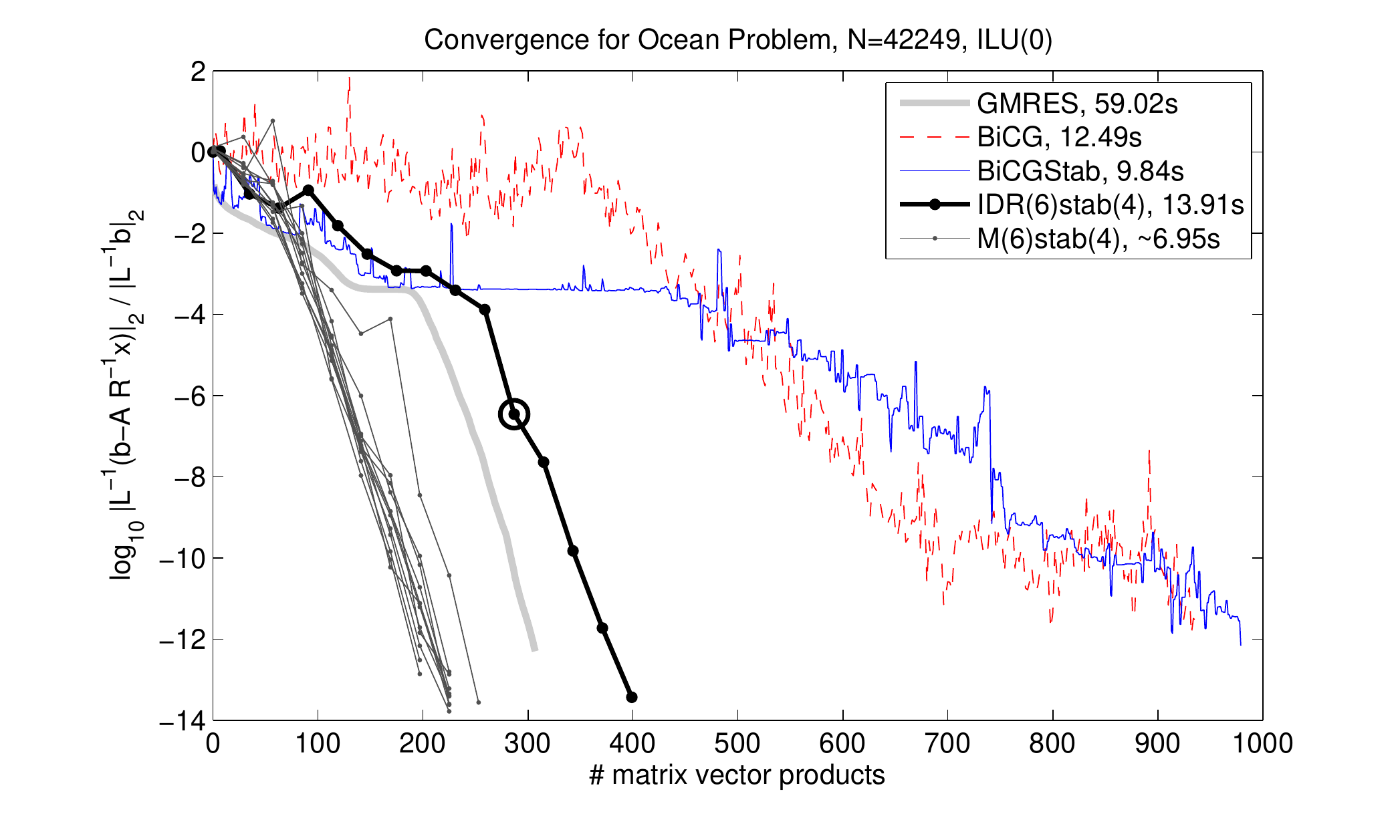}
\caption{\textit{Convergence of \IDRstab, \Mstab and \GMRES for the ocean problem.}}
\label{fig:Ocean2}
\end{figure}

To have an estimate how hard this problem is, the first \rhs is solved with \GMRES, \BiCG and \BiCGstab. Fig.~\ref{fig:Ocean2} shows the convergence of each respective method and gives the computation time in seconds in the legend. Additionally \I solve the preconditioned system with the \IDRstab implementation from Algo. \ref{algo:IDRstab} for $s=6$, $\ell=4$, where after line \ref{algo1:Z24} \I replace the residual by $\br:= \bL^{-1} \cdot (\bb-\bA \cdot \bR^{-1} \cdot \bx)$.

\gI stress the following: 
\begin{enumerate}
	\item In practise one would not solve these systems iteratively, they are by far too small. This is only for test purposes.
	\item \GMRES consumes considerably much time due to long recursions. Anyway, in practise it would not be applicable like that due to storage limits.
	\item For practical problems of large size the condition numbers grow and the preconditioners must be stronger. Then the number of matrix-vector-products dominates the computation time and \IDRstab would out-perform \BiCGstab.
\end{enumerate}

Now we investigate the convergence of \Mstab. To use \Mstab, \we need recycling data in form of matrices $\bU,\bV$. \gI obtain these matrices from the solution process of \IDRstab at the black encircled fetching point in Fig.~\ref{fig:Ocean2}. The obtained $\bU,\bV$ are then used with $\cP$ of \IDRstab as input arguments for \Mstab to solve all twelve \rhs-es subsequently.

From Fig.~\ref{fig:Ocean2} \we see that \Mstab convergences within $200\sim 250$ iterations for each system, whereas \GMRES needs $\approx 300$ and \IDRstab $\approx 400$ iterations. \Mstab solves each system in an average time of $6.95\textrm{s}$, thus achieves a speed-up of $2$ relative to its according \IDRstab variant.

\subsection{Finite Termination of \Mstab}
After \we have seen in the last subsection {that} the residual of \Mstab can drop considerably earlier than for \IDRstab, \I now investigate {why} it does. For this purpose \I present test cases apart from practical applications that show the finite termination behaviour of \Mstab.

Throughout this subsection \we consider the matrix $\bA = \operatorname{tridiag}(2,3,1) \in \R^{N \times N}$, $N=40$ with two \rhs-es $\bb\ha = \boldsymbol{1} \in \R^{N}$, $\bb\hb = \sin(2\cdot \pi/N \cdot (1,...,N)^T\,)\in \R^{N}$. It holds $\bb\ha \perp \bb\hb$.

\begin{figure}
	\centering
	\includegraphics[width=0.75\linewidth]{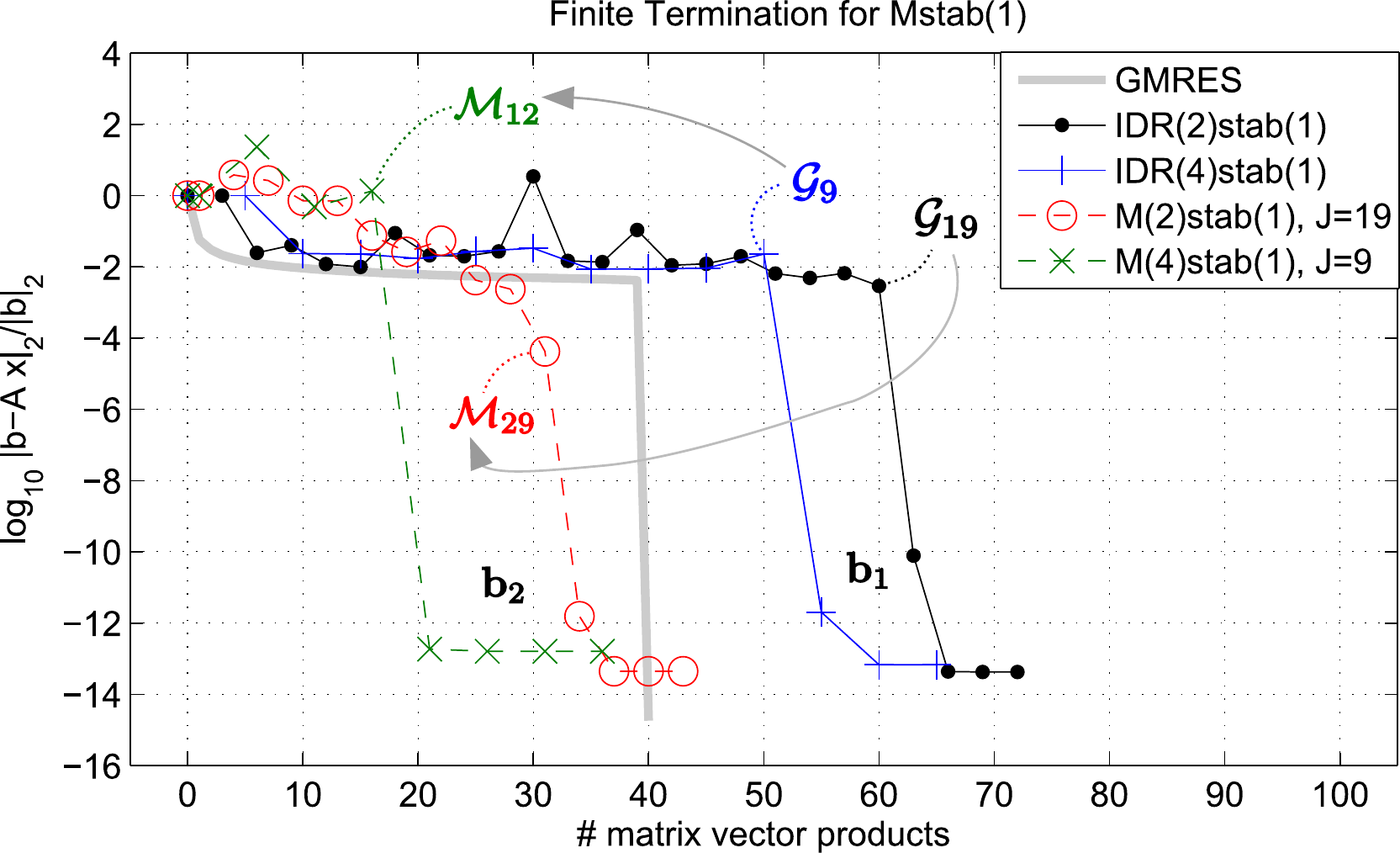}
	\caption{\textit{Termination of \Mstab with $\ell = 1$.}}
	\label{fig:M_Finite_Termination}
\end{figure}

Figs. \ref{fig:M_Finite_Termination} and \ref{fig:Mstab_Termination} show the convergence of \GMRES and different \IDRstab variants for the solution of $\bb\ha$. Additionally, both figures show convergence curves of \Mstab variants.

\paragraph{Termination for $\ell=1$}
Let us first consider in Fig.~\ref{fig:M_Finite_Termination} the solution of $\bb\ha$ with \textsf{IDR($2$)stab($1$)}. Notice that \textsf{IDR($2$)stab($1$)} shoud terminate after a residual in $\cG_{19}$ is computed. This is because of $\cG_{j} = \lbrace \bO \rbrace$ $\forall j > 19$ in exact arithmetic\footnote{assuming the canonical case}, cf. \cite[p. 1050]{IDR-report}.

From the run of \textsf{IDR($2$)stab($1$)} \I fetch recycling data for \textsf{$\cM$($2$)stab($1$)} after $60$ \MV-s (i.e. such that $\opImage(\bV)\subseteq \cG_{19}$). With that data \I then use \textsf{$\cM$($2$)stab($1$)} to solve for $\bb\hb$. From Fig.~\ref{fig:M_Finite_Termination} we see that \textsf{$\cM$($2$)stab($1$)} terminates considerably earlier than \textsf{IDR($2$)stab($1$)}. Assuming the canonical case, the reason for the earlier termination is as follows: Constructing an $\cM$-space from $\cG_{19}$, of which the recycling data is, we have
\begin{align*}
\dim(\cG_{19}) & \leq \max\lbrace 0,\,N - j \cdot s \rbrace = 40 - 19 \cdot 2 = 2\\
\dim(\cM_{19}) & \leq \dim(\cG_{19}) + j \leq 2 + 19 = 21\\
\dim(\cM_{j})  & \leq \max\lbrace 0,\,\dim(\cM_{19}) - (j-19) \cdot 2)\rbrace\,.
\end{align*}
In consequence $\cM_{29}$ is the latest $\cM$-space that differs from $\lbrace \bO \rbrace$. Thus \textsf{$\cM$($2$)stab($1$)} terminates after a residual in $\cM_{29}$ is found.

In Fig.~\ref{fig:M_Finite_Termination} \I show in the same way the convergence of \textsf{IDR($4$)stab($1$)} and \textsf{$\cM$($4$)stab($1$)}. The recycling data was fetched from $\cG_9$. Assuming the canonical case, it holds
\begin{align*}
\dim(\cG_{9}) & \leq 40 - 9 \cdot 4 = 4\\
\dim(\cM_{9}) & \leq \dim(\cG_{9}) + j \leq 4 + 9 = 13\\
\dim(\cM_{j})  & \leq \max\lbrace 0,\,\dim(\cM_{9}) - (j-9) \cdot 4\rbrace\,,
\end{align*}
thus $\cM_{12}$ is the latest non-zero $\cM$-space. In consequence \textsf{IDR($4$)stab($1$)} terminates after the residual is shrinked into $\cM_{12}$.

\paragraph{Termination for general $\ell$}
In Fig.~\ref{fig:Mstab_Termination} \I show cases for different values of $\ell$. First of all \I solve for $\bb\ha$ with \textsf{IDR($4$)stab($1$)} and \textsf{IDR($4$)stab($2$)}. As both methods use a different $\ell$, they work on different Sonneveld spaces. However, we can fetch recycling data from the respective space $\cG_8$ (the indication shows to the fetching points) for both methods. The recycling data from \textsf{IDR($4$)stab($1$)} and from \textsf{IDR($4$)stab($2$)} are denoted by \emph{Data A} and \emph{Data B}, respectively (cf. to the legend).

As both recycling data fullfill the requirements (cf. Algo. \ref{algo:IDRstab}, l.~2) of \textsf{$\cM$($4$)stab($\ell$)}, both data can be used in \textsf{$\cM$($4$)stab($\ell$)} for arbitrary $\ell$. \gI test this as follows: With each recycling data \emph{A} and \emph{B} \I call \textsf{$\cM$($4$)stab($\ell$)} with $\ell = 1$ and $\ell = 2$ to solve for $\bb\hb$. This makes four experiments with four convergence curves, cf. to the legend of Fig.~\ref{fig:Mstab_Termination}.

Assuming the canonical case, we have
\begin{align*}
\dim(\cG_{8}) & \leq 40 - 8 \cdot 4 = 8\\
\dim(\cM_{8}) & \leq \dim(\cG_{8}) + j \leq 8 + 8 = 16\\
\dim(\cM_{j})  & \leq \max\lbrace 0,\,\dim(\cM_{8}) - (j-8) \cdot 4\rbrace\,,
\end{align*}
irrespective of $\ell$, thus $\cM_{11}$ is the latest non-zero $\cM$-space. Indeed, all four curves drop sharply after a residual in $\cM_{11}$ is computed. The numerical results mesh well with the above theory on $\cM$-spaces.

\section*{How to choose the Fetching Point}
\begin{figure}
	\centering
	\includegraphics[width=0.75\linewidth]{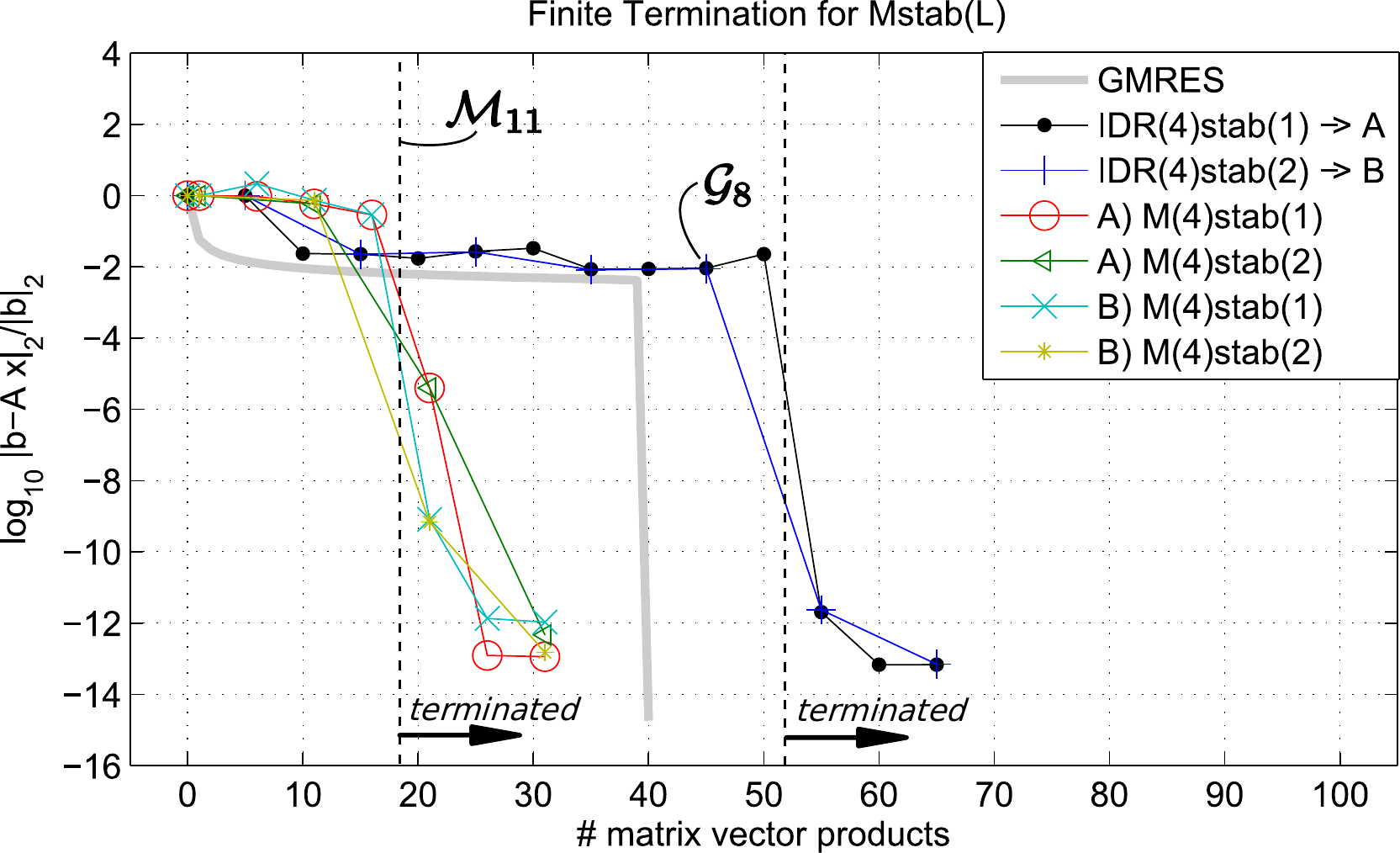}
	\caption{\textit{Termination behaviour of \Mstab for $s=2$ and $\ell \in \lbrace 1,2\rbrace$.}}
	\label{fig:Mstab_Termination}
\end{figure}

In the above finite termination experiments (cf. Figs.~\ref{fig:M_Finite_Termination}-\ref{fig:Mstab_Termination}) the fetching points are chosen such that  the recycling data does not lie in an empty space. This is important. Imagine, e.g., in Fig.~\ref{fig:Mstab_Termination} the recycling data would have been fetched from $\cG_{10}= \lbrace \bO \rbrace$. Then $\bV=\bO$ would hold in theory due to $\opImage(\bV) \subset \cG_j = \lbrace \bO \rbrace$, i.e. in practise we would only recycle round-off. This is not recommended. Instead, in order to achieve a small-dimensional initial $\cM$-space, \I advise to choose the fetching point shortly before the residual drops. 

Finite termination and iterative convergence are quite comparable: E.g., it may be that the full Krylov subspace for $\bb\ha$ of the ocean problem has only $350$ dimensions in exact arithmetic. One even observes that the termination and convergence behaviour of \Mstab, \GMRES and \IDRstab in Fig.~\ref{fig:Mstab_Termination} and \ref{fig:Ocean2} look similar. Thus, \I also recommend for large dimensional systems to choose the fetching point shortly before the residual drops.

In order to estimate a good fetching point in practice, the data may be fetched when half the tolerance is reached. Optionally one can solve two systems with \IDRstab: the convergence of the first system gives an estimate for the convergence of the second system can help in this way, to find a suitable fetching point.

\section{Conclusion}
In this paper \I summarized the theory and algorithm of \IDRstab in an intuitive and graphical way. As side contributions, new increment formulas for Sonneveld spaces and their orthogonal complements were given.

\gI presented a generalization of Sonneveld spaces and the induced dimension reduction theorem that allows for recycling of orthogonality information of test spaces in a simple way. From that \I derived the method \Mstab as a generalization of \IDRstab and Miltenberger's \IDR variant. Numerical experiments confirmed the theoretical termination properties and demonstrated the efficiency of \Mstab (e.g. for the ocean problem the computational effort was halved).

\section{Acknowledgments}
Our special thanks go to Martin van Gijzen for helpful comments in form and content, and for pointing out references \cite{Simoncini}, \cite{Sleijpen1} and \cite{Aihara2}.

\FloatBarrier

\end{document}